\documentclass{amsart}
\usepackage{amsfonts,amsmath,amsthm,amssymb}
\newtheorem{theorem}{Theorem}
\newtheorem{lemma}{Lemma}

\newtheorem{corollary}{Corollary}

\DeclareMathOperator{\charac}{char}
\DeclareMathOperator{\expon}{exp}
\def\gp#1{\langle #1 \rangle}

\title[Lie nilpotency indices of modular  group algebras]
      {Lie nilpotency indices of modular  group algebras}
\author{V.~Bovdi and  J.~B.~Srivastava}


\address{V.~Bovdi,
Institute of Mathematics, University of Debrecen, H-4010 Debrecen,
P.O.Box 12, Hungary
\newline
Institute of Mathematics and Informatics,\quad College of Ny\'\i
regyh\'aza, S\'ost\'oi \'ut 31/b, H-4410 Ny\'\i regyh\'aza, Hungary}
\email{vbovdi@math.klte.hu}
\bigskip
\address{J.~B.~Srivastava\\
Department of Mathematics, Indian Institute of Technology, Delhi,
New Delhi-110016, India } \email{jbsrivas@gmail.com}

\thanks{The research was supported by OTKA grants No. T 43034,
No. K 61007}
 \begin{document}

\begin{abstract}
Let $K$ be a field of positive characteristic $p$ and $KG$ the group
algebra of a group $G$. It is known that if $KG$ is Lie nilpotent
then its upper (or lower) Lie nilpotency index is at most $|G'|+1$,
where $|G'|$ is the order of the commutator subgroup.  The class of
groups $G$ for which these indices are  maximal  or  almost maximal
have already been determined. Here we determine $G$ for which upper
(or lower) Lie nilpotency index is the next highest possible.
\end{abstract}
\maketitle

\section{Introduction and results}
Let $R$ be an associative algebra with identity.  The algebra $R$
can be regarded as a Lie  algebra, called the associated Lie algebra
of $R$, via the  Lie  commutator $[x,y]=xy-yx$, for every $x,y\in
R$. Set $[x_1,\ldots,x_{n}]=[[x_1,\ldots,x_{n-1}],x_{n}],$ where
$x_1,\ldots,x_{n}\in R$. The \emph{$n$-th lower Lie power} $R^{[n]}$
of $R$ is the associative ideal generated by all the Lie commutators
$[x_1,\ldots,x_n]$, where  $R^{[1]}=R$ and $x_1,\ldots,x_n \in R$.
By induction, we define the \emph{$n$-th upper  Lie power} $R^{(n)}$
of $R$ as the associative ideal generated by all the Lie commutators
$[x,y]$, where  $x\in R^{(n-1)}$, $y\in R$ and $R^{(1)}=R$.

An algebra  $R$ is called  \emph{lower Lie nilpotent} (respectively
\emph{upper Lie nilpotent}) if there exists $m$ such that
$R^{[m]}=0$\quad  ($R^{(m)}=0$). The minimal integers $m,n$ such
that $R^{[m]}=0$ and $R^{(n)}=0$  are called \emph{ the lower  Lie
nilpotency index} and \emph{ the upper Lie nilpotency  index} of $R$
and these are denoted by $t_{L}(R)$ and $t^{L}(R)$, respectively.

We would like to mention that  N.~Gupta and F.~Levin \cite{GL} have
given an example of algebra $R$ for which $t_L(R)=3$ but $R$ is not
upper Lie nilpotent. Moreover in that paper, it was shown that if
$R$ is a lower Lie nilpotent ring of class $n$, then the unit group
$U(R)$ is nilpotent of class at most $n$.

An algebra  $R$ is called  \emph{Lie hypercentral} if for each
sequence $\{a_i\}$ of elements of $R$ there exists some $n$ such
that $[a_1,\ldots,a_n]=0$.

Let $U(KG)$ be the  group of units  of a group algebra $KG$. For the
noncommutative modular group algebra $KG$  the following Theorem due
to A.A.~Bovdi, I.I.~Khripta, I.B.S.~Passi, D.S.~Passman and etc.
(see \cite{BKh,K, PPS}) is well known: The following statements are
equivalent: (a) $KG$ is lower Lie nilpotent (b) $KG$ is Lie
hypercentral; (c) $KG$ is upper Lie nilpotent (d) $U(KG)$ is
nilpotent; (e) $\charac(K)=p>0$, $G$ is nilpotent and its commutator
subgroup $G^{\prime}$ is a finite $p$-group.

It is well known (see \cite{SVB}) that,  if $KG$ is Lie nilpotent,
then
$$
t_{L}(KG)\leq t^{L}(KG)\leq \vert G^{\prime}\vert +1.
$$
According to \cite{BP}, if $\charac(K)>3$, then $t_{L}(KG)=
t^{L}(KG)$. But the question of when $t_{L}(KG)= t^{L}(KG)$ for
$\charac(K)=2,3$ is in general still open.

A.Shalev in \cite{S2} began the study of the question when the Lie
nilpotent group algebra $KG$ has  the maximal lower Lie nilpotency
index. In \cite{BS,S2} was given the complete description of such
Lie nilpotent group algebras. In \cite{B1, BJS} we obtained the full
description of the Lie nilpotent group algebras $KG$ with
upper/lower almost maximal Lie nilpotent indices.

Our main results in this paper are  the following theorems.
\begin{theorem}\label{Th1}
Let $KG$ be a Lie nilpotent group algebra  over a field $K$ of
positive characteristic $p$. Then\; $t^L(KG)=\vert
G^{\prime}\vert-4p+5$ \; if and only if one of the following
conditions holds:
\begin{enumerate}
\item[(i)] $p=2$,\quad  $cl(G)=2$\; and   $G^{\prime}\cong C_2\times C_2\times C_2$;
\item[(ii)] $p=5$,\quad  $cl(G)=2$\; and $G^{\prime}\cong C_5\times C_5$.
\end{enumerate}
Moreover,  $t_L(KG)=t^L(KG)$.
\end{theorem}
\begin{theorem}\label{Th2}
Let $KG$ be a Lie nilpotent group algebra  over a field $K$ of
positive characteristic $p$. Then\; $t^L(KG)=\vert
G^{\prime}\vert-3p+4$ \; if and only if one of the following
conditions holds:
\begin{enumerate}
\item[(i)] $p=2$, \quad $cl(G)=3$,\quad  $G^{\prime}\cong
C_2\times C_2\times C_2$ and $\gamma_3(G)$ is cyclic; \item[(ii)]
$p=2$, \quad $G^{\prime}\cong C_4\times C_2$ and
$\gamma_2(G)^2\subseteq \gamma_3(G)$; \item[(iii)] $p=5$,\quad
$cl(G)=3$ and $G^{\prime}\cong C_5\times C_5$.
\end{enumerate}
Moreover,  $t_L(KG)=t^L(KG)$.
\end{theorem}

\begin{theorem}\label{Th3}
Let $KG$ be a Lie nilpotent group algebra  over a field $K$ of
positive characteristic $p$. Then\; $t^L(KG)=\vert
G^{\prime}\vert-2p+3$ \; if and only if one of the following
conditions holds:
\begin{enumerate}
\item[(i)] $p=2$, \quad  $cl(G)=3$,\quad  $G^{\prime}\cong
C_2\times C_2\times C_2$\quad and\quad  $\gamma_3(G)\cong
C_2\times C_2$; \item[(ii)] $p=2$,\quad  $G^{\prime}=\gp{a}\times
\gp{b}\cong C_4\times C_2$ and $\gamma_3(G)$ is one of the
following groups: $\gp{b}$,\quad  $\gp{a^2}\times \gp{b}$,\quad
$\gp{a^2b}$; \item[(iii)] $p=3$,\quad $cl(G)=3$ and
$G^{\prime}\cong C_3\times C_3$.
\end{enumerate}
Moreover,  $t_L(KG)=t^L(KG)$.
\end{theorem}

\begin{theorem}\label{Th4}
Let $KG$ be a Lie nilpotent group algebra  over a field $K$ of
positive characteristic $p$. Then
\begin{itemize}
\item[(i)]\quad $t_L(KG)=\vert G^{\prime}\vert-4p+5$ \quad if and
only if $G$ and $K$ satisfies one of the conditions  (i)--(ii) of
Theorem \ref{Th1}. \item[(ii)] \quad $t_L(KG)= \vert
G^{\prime}\vert-3p+4$ \quad if and only if $G$ and $K$ satisfies one
of the conditions  (i)--(iii) of Theorem \ref{Th2}.
\item[(iii)] \quad $t_L(KG)=\vert G^{\prime}\vert-2p+3$ \quad if
and only if $G$ and $K$ satisfies one of the conditions (i)--(iii)
of Theorem \ref{Th3}.
\end{itemize}
\end{theorem}
Finally, by  Xiankun Du's and  I.~Khripta's Theorems (\cite{D,K}) we
get
\begin{corollary}\label{c:1}
Let $KG$ be the group algebra  of a finite $p$-group $G$ over a
field $K$ of $\charac(K)=p>0$ and $U(KG)$ its group of units. Then
the nilpotency class of  $U(KG)$ is equal to
\begin{itemize}
\item[(i)]\quad $\vert G^{\prime}\vert-4p+4$ \quad if and only if
$G$ and $K$ satisfies one of the conditions  (i)--(ii) of Theorem
\ref{Th1}. \item[(ii)] \quad $\vert G^{\prime}\vert-3p+3$ \quad if
and only if $G$ and $K$ satisfies one of the conditions (i)--(iii)
of Theorem \ref{Th2}. \item[(iii)] \quad $\vert
G^{\prime}\vert-2p+2$ \quad if and only if $G$ and $K$ satisfies one
of the conditions (i)--(iii) of Theorem \ref{Th3}.
\end{itemize}
\end{corollary}

\section{Preliminaries}
We use the  standard   notation:\quad $C_n$ is the cyclic group of
order $n$;\quad  $\zeta(G)$ is the center of a group $G$ \; and \;
$(g,h)=g^{-1}h^{-1}gh$\; ($g,h\in G$); \; $\gamma_{i}(G)$\quad is
the $i$-{th} term of the lower central series of $G$, i.e.
$$
\gamma_{1}(G)=G,\quad \quad
\gamma_{i+1}(G)=\big(\gamma_{i}(G),G\big)\quad \quad (i\geq 1).
$$
Let $K$ be a field of positive  characteristic $p$,\quad  $G$ a
group and put
$$
\mathfrak D_{(m)}(G)=G\cap ( 1+KG^{(m)}),\qquad (m\geq 1).
$$
The subgroup $\mathfrak D_{(m)}(G)$ is  called the $m$-th   \emph{
Lie dimension subgroup} of $KG$ and by Theorem 2.8 (\cite{P}, p.48)
we have:
\begin{equation}\label{e:1}
\mathfrak D_{(m+1)}(G)=\left \{
\begin{array}{ll}
G & \quad \text{if} \quad  m=0;\\
G^{\prime}& \quad \text{if}  \quad     m=1;\\
{\big(\mathfrak D_{(m)}(G),G\big)(\mathfrak D_{(\lceil
{\frac{m}{p}}\rceil +1)}(G))^p}& \quad \text{if}  \quad  m\geq 2,
\end{array} \right.
\end{equation}
where $\lceil {\frac{m}{p}}\rceil $ is the upper integer part of
${\frac{m}{p}}$.

By \cite{P} (see p.46) there exists an explicit expression for
$\mathfrak D_{(m+1)}(G)$:
\begin{equation}\label{e:2}
\mathfrak D_{(m+1)}(G)=\prod_{(j-1)p^i\geq m}\gamma_j(G)^{p^i}.
\end{equation}

Put $p^{d_{(k)}}=[\mathfrak D_{(k)}(G):\mathfrak D_{(k+1)}(G)]$,
where $k\geq 1$. If $KG$ is Lie nilpotent, such that
$|\gamma_2(G)|=p^n$, then according to Jennings' theory \cite{S3}
for the Lie dimension subgroups, we get
\begin{equation}\label{e:3}
\displaystyle t^{L}(KG)=2+(p-1) \sum_{m\geq 1} md_{(m+1)},
\end{equation}
and it is easy to check that
\begin{equation}\label{e:4}
\sum_{m\geq 2} d_{(m)}=n.
\end{equation}
For $x,y,z\in G$  we shall use the following well known formula
\begin{equation}\label{e:5}
(x\cdot y,z)=(x,z)(x,z,y) (y,z).
\end{equation}

We begin with the  following results by A.~Shalev (see Corollary 4.5
and Corollary 4.6 of \cite{S1} and Theorem 3.9 of \cite{S2}):
\begin{lemma}\label{sha} Let $K$ be a field with
$\charac(K)=p>0$ and $G$  a nilpotent group such that
$|G^{\prime}|=p^n$ and  $\expon(G^{\prime})=p^l$.
\begin{enumerate}
\item[(i)]  If $d_{(m+1)}=0$ and $m$ is a power of $p$, then
$\mathfrak D_{(m+1)}(G)=\gp{1}$.
\item[(ii)] If $d_{(m+1)}=0$ and $p^{l-1}$ divides
$m$, then $\mathfrak D_{(m+1)}(G)=\gp{1}$.
\item[(iii)] If $p\geq 5$ and $t_L(KG)<p^n+1$ then $t_L(KG)\leq
p^{n-1}+2p-1$.
\item [(iv)] If $d_{(l+1)}=0$ for some $l<pm$ then  $d_{(pm+1)}\leq d_{(m+1)}$.
\item [(v)] If $d_{(m+1)}=0$ then $d_{(l+1)}=0$ for all $l\geq m$
such  that $\nu_{p'}(l)\geq \nu_{p'}(m)$, where $\nu_{p'}(x)$ is the
maximal divisor of $x$ which relative prime to $p$.
\end{enumerate}
\end{lemma}

Fist of all we proof the following:
\begin{lemma}\label{L:2}
Let $K$ be a field with $\charac(K)=p>0$ and $G$ a nilpotent group
such that  $\vert G^{\prime} \vert=p^{n}$. Then \quad
$t^{L}(KG)=\vert G^{\prime} \vert -4p+5$\quad if and only if one of
the following condition holds:
\begin{itemize}
\item [(i)]  $p=2$, \quad $n=3$\quad  and \quad $d_{(2)}=3$; \item
[(ii)] $p=3$, \quad  $n=3$ \quad  and \quad
$d_{(2)}=d_{(4)}=d_{(6)}=1$; \item [(iii)] $p=5$, \quad  $n=2$
\quad and \quad  $d_{(2)}=2$.
\end{itemize}
\end{lemma}

\begin{proof}
Let $t^{L}(KG)=\vert G^{\prime} \vert -4p+5$, where $p=\charac(K)$.
If either $n=1$ or $p=n=2$,  then  according to Theorem 1 of
\cite{BS}, and to Corollary  1 of \cite{BJS}, we get $t^L(KG)\geq
\vert G^{\prime} \vert$. So in the sequel we assume that either
$p=2$ and $n\geq 3$ or  $n\geq 2$.

Let $p=2$, $n=3$ and suppose that $d_{(2)}\leq 2$. By (\ref{e:2}) we
have that $t^L(KG)\geq 6> 5$, so only possible case  is $d_{(2)}=3$.

Let $p=2$, $n=4$ and consider the following cases: $d_{(2)},
d_{(3)}\in \{ 1,2,3\}$. By (\ref{e:2}), (\ref{e:3}) and (i) of Lemma
\ref{sha} it is easy to compute that one  possible solution is:
\quad $d_{(2)}=d_{(3)}=1$\quad  and\quad $d_{(5)}=2$ which in
contradictions  with (iv) of Lemma \ref{sha}.

Now, let $p=2$ and $n\geq 5$.  We prove that $d_{(p^{i}+1)}>0$
for\quad $0\leq i\leq n-2$.

Suppose  that  $d_{(2^{n-2}+1)}=0$ and $d_{(2^{n-3}+1)}\not=0$. By
(i) of Lemma \ref{sha} we have that $\mathfrak
D_{(2^{n-2}+1)}(G)=\gp{1}$ and $d_{(r)}=0$ for every $r\geq
2^{n-2}+1$. Moreover, if $d_{({q}+1)}\not=0$, then $q<2^{n-2}$.
According to (\ref{e:3}) it follows that
\begin{align*}
t^{L}(KG)&=2+\sum_{i=0}^{n-3}{2^{i}}+\sum_{i=0}^{n-3}{2^{i}}(d_{({2^{i}}+1)}-1)+ \sum_{q\not=2^i}qd_{({q}+1)}\\
&<1+(1+\sum_{i=0}^{n-3}(d_{({2^{i}}+1)}-1)+ \sum_{q\not=2^i}d_{({q}+1)})\cdot 2^{n-2}\\
&=1+(1+n-(n-2))\cdot 2^{n-2} < 2^n-3,
\end{align*}
which is  contradicts to $t^{L}(KG)=\vert G^{\prime} \vert -3$.

Therefore \; $d_{(p^{i}+1)}>0$\; for\; $0\leq i\leq n-2$ and  by
(\ref{e:3}) and (\ref{e:4}) there exists $\alpha\geq 2$ such that
$d_{(\alpha+1)}=1$, and
\begin{align*}
t^{L}(KG)= 2+\sum_{i=0}^{n-2}{2^{i}}+\alpha d_{(\alpha+1)}=
1+2^{n-1}+\alpha.
\end{align*}
Since $t^{L}(KG)=|G'|-3$, it must be $\alpha =2^{n-1}-4$. Put
$m=2^{n-3}-1$ and $l=4m$. Since $\nu_{2'}(l)=\nu_{2'}(m)$ and
$d_{((2^{n-3}-1)+1)}=0$, by (v) of Lemma \ref{sha} we get
$d_{((2^{n}-4)+1)}=0$, a contradiction.

Let now $p=3$ and $n\geq 4$. Using  the same arguments as in the
previous case, we get $d_{(3^{i}+1)}=d_{((3^{n-1}-4)+1)}=1$ for
$0\leq i\leq n-2$. Since \quad $exp(G^{\prime})\leq 3^{n-1}$ \quad
and \quad $3^{n-2}\;\mid\;2\cdot 3^{n-2}$,\quad  by (ii) of Lemma
\ref{sha} it yields  that $\mathfrak D_{(2\cdot
3^{n-2}+1)}(G)=\gp{1}$, so $d_{((3^{n-1}-4)+1)}=0$, a contradiction.

If $n=2$, then by (\ref{e:3}) we obtain that $t^L(KG)\geq 6>2$.

Let $p=n=3$. According to the case when $p=2$ and $n=4$, it is easy
to check that $d_{(2)}=d_{(4)}=d_{(6)}=1$ and the proof for $p=3$ is
complete.

Let  $p=5$. The inequality $t^L(KG)=5^n-15\leq 5^{n-1}+9$ is
verified for $n=2$. Suppose  $d_{(2)}=1$. According to (\ref{e:3})
we get  $ t^L(KG)\geq14 >10$.

Finally, let $p\geq 7$. Since $p^{n-1}>6$,\quad
$p^{n}-4p+5>p^{n-1}+2p-1$, so by (iii) of Lemma \ref{sha} it follows
that  this case is impossible. \end{proof}

\begin{lemma}\label{L:3}
Let  $K$ be a field with  $\charac(K)=p\geq 3$ and $G$  a
nilpotent group, such that $\vert G^{\prime}\vert=p^n$  and
$t^{L}(KG)=p^n-4p+5$. Then $p=5$, $cl(G)=2$ and $G'\cong C_5\times
C_5$.
\end{lemma}
\begin{proof} Let $\charac(K)=p\geq 3$ and $t^{L}(KG)=p^n-4p+5$. Clearly, either
(iii) or  (iv) of Lemma \ref{L:2} holds, so  we consider the
following cases:

Case {\bf 1}. Let $p=n=3$ and   $d_{(2)}=d_{(4)}=d_{(6)}=1$. If
$cl(G)=2$ by Theorem 3.2 of \cite{BK} we have that
$t^L(KG)=t(G^{\prime})+1\leq 12<20$. So assume that $cl(G)\geq 3$.
If $exp(G^{\prime})=3$ we obtain that $\mathfrak
D_{(3)}(G)=\gamma_3(G)$ and $\mathfrak D_{(4)}(G)=\gamma_4(G)$, so
$\gamma_3(G)=\gp{1}$, a contradiction.

Now let $exp(G^{\prime})=3^2$. If $cl(G)=3$ then $\mathfrak
D_{(2)}(G)=\gamma_2(G)\cong C_9\times C_3$ and $\mathfrak
D_{(3)}(G)=\gamma_3(G)\cdot \gamma_2(G)^3=\{C_9, C_3 \times C_3\}$.

If $\mathfrak D_{(3)}(G)$  is not cyclic we have that $\mathfrak
D_{(4)}(G)=(\mathfrak D_{(3)}(G),G)=\mathfrak D_{(3)}(G)$ which is
impossible. In the other case, we obtain that
$$
\mathfrak D_{(4)}(G)=(\mathfrak D_{(3)}(G),G)\cdot \mathfrak
D_{(3)}(G)^3=\mathfrak D_{(3)}(G)
$$
and again it holds that  $(\mathfrak D_{(3)}(G),G)=\mathfrak
D_{(3)}(G)$.

Finally suppose that $cl(G)=4$. Clearly $\mathfrak
D_{(3)}(G)=\gamma_3(G)$ and
$$
\mathfrak D_{(4)}(G)=(\mathfrak D_{(3)}(G),G)\cdot
\gamma_2(G)^3=\gamma_4(G)\cdot\gamma_2(G)^3=\gamma_3(G),
$$
so $\gamma_2(G)^3\subseteq \gamma_3(G)$ and by Theorem III.2.13
(\cite{H}, p.266), we have that $\gamma_3(G)^3 \subseteq
\gamma_{4}(G)$. Now $ \mathfrak D_{(5)}(G)=\gamma_4(G)\cdot
\gamma_3(G)^3=\gamma_4(G)$ and
$$
\mathfrak D_{(6)}(G)=\gamma_5(G)\cdot \gamma_3(G)^3=\gamma_4(G),
$$
so $\gamma_3(G)\cong C_9$  and $\mathfrak D_{(7)}(G)=
\gamma_3(G)^3=\mathfrak D_{(6)}(G)$, which is a contradiction.

Case {\bf 2}. Let $p=5$,  $n=2$ and  $d_{(2)}=2$. Obviously,
$\gamma_2(G)\cong C_5 \times C_5$. If $cl(G)=3$, then $\mathfrak
D_{(3)}(G)=\gamma_3(G)\not=\gp{1}$, a contradiction.\end{proof}

{\bf Proof of the Theorem \ref{Th1}}.  Follows from   Lemmas
\ref{L:2}- \ref{L:3}. The equality $t_L(KG)=t^L(KG)$ is a
consequence of (ii) of Theorem 3.2 of \cite{BK}.
\bigskip

\begin{lemma}\label{L:4}
Let $K$ be a field with $\charac(K)=p>0$ and $G$ a nilpotent group
such that  $\vert G^{\prime} \vert=p^{n}$. Then \quad
$t^{L}(KG)=\vert G^{\prime} \vert -3p+4$\quad if and only if one of
the following condition holds:
\begin{itemize}
\item [(i)]  $p=2$, \quad $n=3$,\quad   $d_{(2)}=2$\quad and \quad
$d_{(3)}=1$; \item [(ii)] $p=3$, \quad  $n=3$ \quad and \quad
$d_{(2)}=d_{(4)}=d_{(7)}=1$; \item [(iii)] $p=5$, \quad  $n=2$
\quad and \quad  $d_{(2)}=d_{(3)}=1$.
\end{itemize}
\end{lemma}
\begin{proof}
As in Lemma \ref{L:2} we can assume that either $p=2$ and $n\geq 3$
or $n\geq 2$.

Let $p=2$, $n=3$ and suppose that $d_{(2)}=1$. By (\ref{e:3}) we
have that $t^L(KG)\geq 7> 6$, so only possible case  is $d_{(2)}=2$
and by (i) of Lemma \ref{sha} the statement is holds.

Let $n\geq 4$. Using the  same arguments as in the proof of the
Lemma \ref{L:2} we obtain that $d_{(2^{i}+1)}=d_{((2^{n-1}-3)+1)}=1$
\; and \; $d_{(j)}=0$\;, where \; $0\leq i\leq n-2$, \; $j\neq 2^i
+1$, \; $j\neq 2^{n-1}-2$\; and \; $j>1$.  The subgroup $H=\mathfrak
D_{( 2^{n-1}-2)}(G)$ is central of order $2$ and   from (\ref{e:2})
it yields
\begin{equation}\label{e:6}
\begin{aligned}
\mathfrak D_{(m+1)}(G)/H &=\prod_{(j-1)2^i\geq
m}\gamma_j(G)^{2^i}/H\\
&=\prod_{(j-1)2^i\geq m}\gamma_j(G/H)^{2^i}=\mathfrak
D_{(m+1)}(G/H).
\end{aligned}
\end{equation}
Put $2^{\overline{d}_{(k)}}=[\mathfrak D_{(k)}(G/H):\mathfrak
D_{(k+1)}(G/H)]$ for $k\geq 1$. It is easy to check that \;
$\overline{d}_{(2^{i}+1)}=1$\; and \; $\overline{d}_{(j)}=0$, where
\; $0\leq i\leq n-2$, \; $j\neq 2^i +1$\; and \; $j>1$.

Clearly, $|\gamma_2(G/H)|=2^{n-1}$ and
$t^L(K[G/H])=|\gamma_2(G/H)|+1$. So by Lemma 3 of \cite{BS} and by
Theorem 1 of \cite{BS} the group $\gamma_2(G/H)$ is either a
cyclic $2$-group or $C_2\times C_2$. If $\gamma_2(G/H)$ is  a
cyclic $2$-group, then  by (a) of Lemma III.7.1 (\cite{H}, p.300)
we have that   $\gamma_2(G)$ is abelian, so it is isomorphic to
either $C_{2^{n-1}}\times C_2$ or $C_{2^{n}}$. If $\gamma_2(G)$ is
cyclic, then by Theorem 1 of \cite{BS} we get
$t^L(KG)=|\gamma_2(G)|+1$, so we do not consider this case. On the
other hand, if $\gamma_2(G/H)\cong C_2\times C_2$, then
$|\gamma_2(G)|=8$. Since there are no nilpotent groups with
nonabelian commutator subgroup of order $8$ (see \cite{BR}), we
can put that $\gamma_2(G)=\gp{a,b\mid a^{2^{n-1}}=b^2=1}\cong
C_{2^{n-1}}\times C_2$.

According to (1) it holds
$$
\mathfrak D_{(2)}(G)=\gamma_2(G),\qquad\qquad \mathfrak
D_{(3)}(G)=\gamma_3(G)\cdot\gp{a^2}.
$$
Since   $\vert \mathfrak D_{(2)}(G)/\mathfrak D_{(3)}(G)\vert =2$ we
obtain one of the following cases:
\[
\begin{aligned}
\gamma_3(G)&=\gp{a},\qquad\quad \gamma_3(G)=\gp{ab},\qquad\qquad\gamma_3(G)=\gp{b},\\
 \gamma_3(G)&=\gp{a^{2^j}b},\qquad \gamma_3(G)=\gp{a^{2^j},b}, \qquad (1\leq j\leq n-2).\\
\end{aligned}
\]
 We consider each of these cases:

Case {\bf 1}.  Let either $\gamma_3(G)=\gp{a}$ or
$\gamma_3(G)=\gp{ab}$. Since $\gamma_2(G)^2\subset \gamma_3(G)$, by
Theorem III.2.13 (\cite{H}, p.266), we have that $\gamma_k(G)^2
\subseteq \gamma_{k+1}(G)$ for every $k\geq 2$. It follows that
$\gamma_2(G)^2=\gamma_4(G)$. Moreover, $\gamma_{3}(G)^2\subseteq
\gamma_{5}(G)$. Indeed, the elements of the form $(x,y)$, where
$x\in \gamma_2(G)$ and $y\in G$ are generators of $\gamma_3(G)$, so
we have to  prove that $(x,y)^2 \in \gamma_5(G)$. By (5)
$$
(x^2,y)=(x,y)(x,y,x)(x,y)=(x,y)^2(x,y,x)^{(x,y)}
$$
and  $(x^2,y),\; (x,y,x)^{(x,y)}\in \gamma_5(G)$, so $(x,y)^2\in
\gamma_5(G)$ and $\gamma_3(G)^2\subseteq\gamma_5(G)$. Thus
$|\gamma_3(G)|=2$,  a contradiction.

Case {\bf 2}. Let $\gamma_3(G)=\gp{b}$. Now, let us compute the weak
complement of $\gamma_3(G)$ in $\gamma_2(G)$ (see \cite{BK}, p.34).
It is easy to see that, with the notation of \cite{BK} (see p.34)
$$
P=\gamma_2(G),\qquad  H=\gp{b},\qquad  H\setminus P^2=\{b\}
$$
so $\nu(b)=2$ and  the weak complement is $A=\gp{a}$. Since $G$ is
of class $3$, by (ii) of Theorem 3.3 (\cite{BK}, p.43) we have
\[
\begin{split}
t_L(KG)=t^L(KG)&=t(\gamma_2(G))+t(\gamma_2(G)/\gp{a})\\
&=2^{n-1}+3\not= |G^\prime|-2.
\end{split}
\]

Case {\bf 3}.  Let $\gamma_3(G)=\gp{a^{2^{j}}b}$ with  $1\leq
j\leq n-2$.  Since $d_{(3)}=1$,  by (\ref{e:1})  we get
$$
\mathfrak D_{(2)}(G)=\gamma_{2}(G), \quad \mathfrak
D_{(3)}(G)=\gp{a^2,b},\quad \mathfrak
D_{(4)}(G)=(\gp{a^2,b},G)\cdot\gp{a^{4}},
$$
and $\mathfrak D_{(4)}(G)$ is one of the following groups:
$\gp{a^{2}}$,\; $\gp{a^{2}b}$,\; $\gp{a^{4},b}$.

Suppose that $\mathfrak D_{(4)}(G)=\gp{a^{2}}$. Since $d_{(4)}=0$ we
have
$$
\mathfrak D_{(5)}(G)=(\gp{a^{2}},G)\cdot \mathfrak
D_{(3)}(G)^2=(\gp{a^{2}},G)\cdot \gp{a^{4}}=\gp{a^{2}}.
$$
The last equality forces $(\gp{a^{2}},G)=\gp{a^{2}}$ and so
$\mathfrak D_{(k)}(G)=\gp{a^{2}}$ for each $k\geq 5$, which is
impossible.

Now let  $\mathfrak D_{(4)}(G)=\gp{a^{2}b}$. Similarly, because
$d_{(4)}=0$, we get
$$
\mathfrak D_{(4)}(G)=\mathfrak D_{(5)}(G)=(\gp{a^{2}b},G)\cdot
\gp{a^{4}}=\gp{a^{2}b}
$$
and  $(\gp{a^{2}b},G)=\gp{a^{2}b}$, which is impossible.

Finally, suppose that   $\mathfrak D_{(4)}(G)=\gp{a^{4}}\times
\gp{b}$ and  exists $k \leq 2^{n-2}+1$, such that $\mathfrak
D_{(k)}(G)$ is cyclic. Using the same arguments of above, we obtain
that $\mathfrak D_{(m)}(G)\not=\gp{1}$ for each $m$, which is
impossible.

Therefore $\mathfrak D_{(2^{n-2}+1)}(G)=\gp{a^{2^{n-2}}}\times
\gp{b}$. Since $d_{(2^{n-2}+1)}=1$  and $d_{(2^{n-2}+2)}=0$,  by
(\ref{e:1}) we get
\[
\begin{aligned}
\mathfrak D_{(2^{n-2}+2)}(G)&=(\mathfrak D_{(2^{n-2}+1)}(G), G
)\cdot
\mathfrak D_{(2^{n-3}+2)}(G)^2\\
&=(\mathfrak D_{(2^{n-2}+1)}(G), G )=\gp{\omega\mid \omega^2=1};\\
\mathfrak D_{(2^{n-2}+3)}(G)&=(\mathfrak
D_{(2^{n-2}+2)}(G),G)=(\gp{\omega},G)=\gp{\omega}.
\end{aligned}
\]
It  follows that $(\gp{\omega},G)=\gp{\omega}$ and $\mathfrak
D_{(k)}(G)=\gp{\omega}$ for all $k\geq 2^{n-2}+2$, again  a
contradiction.

Case {\bf  4}. Let $\gamma_3(G)=\gp{a^{2^{j}}}\times \gp{b}$ with
$1\leq j\leq n-2$. It is easy to check that this case is similar to
the last subcase of the previous case.

Therefore the case  $p=2$ is finished.

Let now $p=3$ and $n\geq 4$. Using  the same arguments as in the
previous case, we get $d_{(3^{i}+1)}=d_{((3^{n-1}-3)+1)}=1$ for
$0\leq i\leq n-2$. Since \quad $exp(G^{\prime})\leq 3^{n-1}$ \quad
and \quad $3^{n-2}\;\mid\;2\cdot 3^{n-2}$,\quad  by (ii) of Lemma
\ref{sha} it yields  that $\mathfrak D_{(2\cdot
3^{n-2}+1)}(G)=\gp{1}$, so $d_{((3^{n-1}-3)+1)}=0$, a contradiction.

Let $p=3$ and $n=3$. Similarly to prof of part (ii) of Lemma
\ref{L:2}, we obtain that $d_{(2)}=d_{(4)}=d_{(7)}=1$, so this
case  is complete.

For $p=5$, the part  (iii) of Lemma \ref{sha} gives  only one
possible case:  $n=2$, which follows that  $d_{(2)}=1$ and
$d_{(3)}=1$.
\end{proof}

\begin{lemma}\label{L:5}
Let  $K$ be a field with  $\charac(K)=2$ and $G$  a nilpotent group
such that $\vert G^{\prime}\vert=2^n$. If $t^{L}(KG)=2^n-2$, then
either (i)  or (ii) of Theorem \ref{Th2} hold.
\end{lemma}
\begin{proof}
By (i) of Lemma  \ref{L:4} we have $n=3$, \quad $d_{(2)}=2$\quad and
$d_{(3)}=1$. If $\gamma_2(G)$ is elementary abelian, according to
Theorem {1} and Theorem {1} of \cite{BJS} we have that $cl(G)=3$. If
$|\gamma_3(G)|=4$, then it is easy to check $d_{(2)}=1$, a
contradiction. Therefore $\gamma_2(G)\cong C_4\times C_2$ and  the
statement follows at once by easy calculation.
\end{proof}

\begin{lemma}\label{L:6}
Let  $K$ be a field with  $\charac(K)=p\geq 3$ and $G$  a nilpotent
group, such that $\vert G^{\prime}\vert=p^n$. If
$t^{L}(KG)=p^n-3p+4$ then  $p=5$, $cl(G)=3$ and $\gamma_2(G)\cong
C_5\times C_5$.
\end{lemma}
\begin{proof} Using the same argument of proof of Lemma \ref{L:3},
we obtain that $p\not=3$. If $p=5$ and $n=2$, then \quad
$d_{(2)}=d_{(3)}=1$\quad  and it follows  that $\gamma_2(G)$ is
not central, so  the proof is complete.\end{proof}

{\bf Proof of the Theorem 2}.  Follows from   Lemma \ref{L:5}, Lemma
\ref{L:6} and Theorem \ref{Th1}. The equality $t_L(KG)=t^L(KG)$ is
an immediate consequence of part (ii) of Theorem 3.2 by \cite{BK}
and part (ii) of Theorem 3.3 by \cite{BK}.
\bigskip

\begin{lemma}\label{L:7}
Let $K$ be a field with $\charac(K)=p>0$ and $G$ a nilpotent group
such that  $\vert G^{\prime} \vert=p^{n}$. Then \quad
$t^{L}(KG)=\vert G^{\prime} \vert -2p+3$\quad if and only if one of
the following condition holds:
\begin{itemize}
\item [(i)]  $p=2$, \quad $n=3$,\quad  $d_{(2)}=1$\quad and \quad
$d_{(3)}=2$; \item [(ii)] $p=3$, \quad  $n=2$ \quad and \quad
$d_{(2)}=2$.
\end{itemize}
\end{lemma}
\begin{proof} Let $t^{L}(KG)=\vert G^{\prime} \vert -2p+3$.
By (iii) of  Lemma \ref{sha} we can assume that  $p\leq 3$.

Let $p=2$ and  $n=3$. According to  (i) of  Lemma \ref{L:2},  (i)
of Lemma \ref{L:4} and (\ref{e:3}) it is  easy to check that
$d_{(2)}=1$ and $d_{(3)}=2$.

If  $n\geq 4$, then using the same technic as in the proof of
Lemma \ref{L:2} and Lemma \ref{L:7} we obtain that \quad
$d_{(2^{i}+1)}=d_{(2^{n-1}-1)}=1$ \quad and \quad $d_{(j)}=0$,
where \quad $0\leq i\leq n-2$, \quad $j\neq 2^i +1$, \quad $j\neq
2^{n-1}-1$\quad and $j>1$;

Put $m=2^{n-2}-1$ and $l=2m$. Since $d_{((2^{n-2}-1)+1)}=0$\quad
and\quad $\nu_{2'}(l)=\nu_{2'}(m)$,\quad  by (v) of Lemma \ref{sha}
we get $d_{((2^{n-1}-2)+1)}=0$, a contradiction.

Finally, let $p=3$ and $n \geq 3$. In this case as in proof of (ii)
of  Lemma \ref{L:2}, we obtain the contradiction and the proof is
complete.
\end{proof}

\begin{lemma}\label{L:8}
Let  $K$ be a field with  $\charac(K)=2$ and $G$  a nilpotent group
such that $\vert G^{\prime}\vert=2^n$. If $t^{L}(KG)=2^n-1$, then
either (i) or (ii) of Theorem \ref{Th3} holds.
\end{lemma}
\begin{proof}
If  $n\geq 4$, then the proof is the same of that used in Lemma
\ref{L:3} for $n\geq 5$. So assume that $n=3$.

Therefore either $\gamma_2(G)\cong C_2\times C_2\times C_2$ or
$\gamma_2(G)\cong C_4\times C_2$. In the first case,  part (i) of
Theorem  \ref{Th1} and part (i) of Theorem
 \ref{Th2}   force that  $\gamma_3(G)$ to be of order $4$
and central. In the last one since\quad  $|\mathfrak
D_{(2)}(G)/\mathfrak D_{(3)}(G)|=2$ \quad  and \quad  $\mathfrak
D_{(3)}(G)=\gamma_3(G)\cdot \gamma_2(G)^2$ we obtain that
$\gamma_{3}(G)$ is one of the following groups:
\[
\begin{aligned}
\gamma_3(G)=\gp{a},&\qquad \gamma_3(G)=\gp{ab},\qquad\gamma_3(G)=\gp{b},\\
\gamma_3(G)&=\gp{a^{2}b},\qquad \gamma_3(G)=\gp{a^{2},b}.\\
\end{aligned}
\]

If $\gamma_3(G)=\gp{a}$ or  $\gamma_3(G)=\gp{ab}$ the proof was
made in Case 1 of Lemma \ref{L:3}.

Let either $\gamma_3(G)=\gp{b}$ or $\gamma_3(G)=\gp{a^{2}b}$. It
is easy to check that in both cases  the week complement of
$\gamma_3(G)$ in $\gamma_2(G)$ (see \cite{BK}, p.34) is $A=\gp{a}$
and $t^L(KG)=t_L(KG)=7$.

Finally, let $\gamma_3(G)$ be noncyclic. Obviously, $\mathfrak
D_{(3)}(G)=\gamma_3(G)$ and $\gamma_3(G)^2=\gp{1}$, so $\mathfrak
D_{(4)}(G)=\gp{1}$ and the statement holds.\end{proof}

{\bf Proof of the Theorem 3}. Let  $K$ be a field with
$\charac(K)=p\geq 3$. By Lemma \ref{L:7} it holds that $\gamma_2(G)$
is noncyclic of order $9$. If $p=2$ the statement was proved in
Lemma \ref{L:8}.

{\bf Proof of the Theorem \ref{Th4}}. Let  $p=2$ and
$t_L(KG)=|G'|-2p+3$. Note that in \cite{B1, BS,BJS} we obtained the
description of Lie nilpotent group algebras $KG$ with\quad
$t^L(KG)\in\{|G'|-p+2,\quad  |G'|+1\}$.\quad  Since $t_L(KG)\leq
t^L(KG)$ and \qquad $|G'|-2p+3<|G'|-p+2< |G'|+1$,\qquad  we have
that
$$
t^L(KG)\in\{\;|G'|-2p+3,\quad  |G'|-p+2,\quad  |G'|+1\;\}.
$$
If $t^L(KG)\in\{|G'|-p+2, \quad |G'|+1\}$, then according to
\cite{B1, BS,BJS}, respectively, we have that \quad  $
t_L(KG)\in\{|G'|-p+2, \quad |G'|+1\}$, a contradiction. Therefore
$t^L(KG)=|G'|-2p+3$ and by Theorem \ref{Th3} we get
$t_L(KG)=t^L(KG)$ and the description of groups $G$ is done.

Let $p=3$. First, by (\ref{e:3}) there is no $KG$ with
$t^L(KG)=|G'|-2s$ with  $s\in \mathbb N$, because
$t^L(KG)=2+2(\sum_{q\geq 1}qd_{q+1})\not=3^{n}-2s$. Using this
remark,  the proof is the same as for case $p=2$.

Now the proof of the parts (i) and (ii) of the Theorem \ref{Th4} for
$p=2$ and $p=3$ is the same. Note that for  $p\geq 5$, the proof is
a trivial consequence of the Bhandari-Passi's Theorem \cite{BP}.


\bibliographystyle{abbrv}
\bibliography{Bovdi_Srivastava_final}

\end{document}